\documentstyle[12pt]{article}
\begin{document}

\title{\bf On sufficient conditions for the total
positivity
and for the multiple positivity of matrices. }
\author{Olga M. Katkova, Anna M. Vishnyakova\\
Dept. of Math., Kharkov State University, \\
Svobody sq., 4, 61077, Kharkov, Ukraine, \\
e-mail: olga.m.katkova@ilt.kharkov.ua , \\
anna.m.vishnyakova@univer.kharkov.ua }
\maketitle

\begin{abstract}
The following theorem is proved.

{\bf Theorem.}  { \it  Suppose $M = (a_{i,j})$ be a  $k \times k$
matrix with positive entries and $a_{i,j}a_{i+1,j+1} >  4\cos ^2
\frac{\pi}{k+1}\ a_{i,j+1}a_{i+1,j} \quad (1 \leq i \leq k-1,\ 1
\leq j \leq k-1).$ Then $\det M > 0 .$}

The constant $4\cos ^2 \frac{\pi}{k+1}$ in this Theorem is sharp.
A few other results concerning totally positive and multiply
positive matrices are obtained.
\end{abstract}


{\it 2000 Mathematics Subject Classification  15A48, 15A57, 15A15.}

\section{Introduction and statement of results.}

This paper is inspired by the interesting work \cite{cc} in
which some useful and easily verified conditions of strict total
positivity of a matrix are obtained. We recall that a matrix $A$
is said to be $k$-times positive, if all minors of $A$ of order
not greater than $k$ are non-negative. A matrix $A$ is said to be
multiply positive if it is $k$-times positive for some $k\in {\bf N}.$
A matrix $A$ is said to be totally positive, if all minors of $A$ are
non-negative. For more information about these notions and their
applications we refer the reader to \cite{ando} and \cite{tp}.
According to \cite{tp}
we will denote the class of all $k$-times positive matrices by
$TP_k$ and the class of all totally positive matrices by $TP.$ By
$STP$ we will denote the class of matrices with all minors being
strictly positive and by $STP_k$ the class of matrices with all
minors of order not greater than $k$ being strictly positive.

In \cite{cc} the following theorem was proved

{\bf Theorem A.} {\it Denote by $\tilde{c}$  the unique real root
of $x^3-5x^2+4x-1 = 0$ ($\tilde{c} \approx 4{.} 0796$).
Let  $M = (a_{i,j})$ be an  $n \times n$ matrix with the property that \\
(a) $ a_{i,j} > 0 \  (1 \leq i,j \leq n)$ and \\
(b) $a_{i,j}a_{i+1,j+1} \geq \tilde{c}\  a_{i,j+1}a_{i+1,j} \
(1 \leq i,j \leq n-1) .$\\
Then $M$ is strictly totally positive. }

Note that the verification of total positivity is, in general, a
very difficult problem. Surely, it is not difficult to calculate the
determinant of a given matrix with numerical entries. But
if the order of a matrix or the entries of a
matrix depend on some parameters then the testing of multiple
positivity is complicated. Theorem A provides a convenient
sufficient condition for total positivity of a matrix.

For $c \geq 1$ we will denote by $TP_2(c)$ the class of all matrices $M =
(a_{i,j})$  with positive entries which satisfy the condition
\begin{equation}
\label{def1} a_{i,j}a_{i+1,j+1} \geq c\  a_{i,j+1}a_{i+1,j} \quad
\mbox{{\it for all}} \quad i,j .
\end{equation}

For $c \geq 1$ we will denote by $STP_2(c)$ the class of all matrices $M =
(a_{i,j})$  with positive entries which satisfy the condition
\begin{equation}
\label{def2} a_{i,j}a_{i+1,j+1} > c\  a_{i,j+1}a_{i+1,j}
\quad \mbox{{\it for all}} \quad i,j .
\end{equation}

It is easy to verify that $STP_2 = STP_2(1).$ Theorem A states
that $TP_2(\tilde{c}) \subset STP .$

Denote by
$$c_k := 4 \cos ^2 \frac
{\pi}{k+1}, \ k=2,3,4, \ldots \ .$$

The main result of this paper is the following:

{\bf Theorem 1.}   {\it  Suppose $M = (a_{i,j})$ be a  $k \times k$
matrix with positive entries.  \\
(i) if $M \in TP_2(c_k)$  then $\det M \geq 0$; \\
(ii) if $M \in STP_2(c_k)$  then $\det M >0 $.}

In the proof of Theorem 1 we will show that if $M\in TP_2(c)$ then
every submatrix of $M$ belongs to $TP_2(c).$ Therefore the
following theorem is the simple consequence of Theorem 1.

{\bf Theorem 2.}   {\it  For every $c \geq c_k $ we have \\
(i) if $M \in TP_2(c)$  then $M \in TP_k $; \\
(ii) if $M \in STP_2(c)$  then $M \in STP_k $.}

The following fact is a simple consequence of this theorem.

{\bf Theorem 3.}   {\it   For every $c \geq 4$ we have \\
if $M \in TP_2(c)$  then $M \in STP .$ }

The following statement demonstrates that the constants in
Theorems 1 and 3 are unimprovable not only in the class of
matrices with positive entries but in the classes of Toeplitz
matrices and of Hankel matrices. We recall that a matrix $M$ is
a Toeplitz matrix if it is of the form $M=(a_{j-i})$ and a matrix
$M$ is a Hankel matrix if it is of the form $M=(a_{j+i}).$

{\bf Theorem 4.}{\it \\
(i) For every $ 1 \leq c < c_k $ there
exists a $k\times k $ Toeplitz matrix  $M \in TP_2(c) $
with $\det M < 0 ;$\\
(ii)  for every $ 1 \leq c < c_k $ there exists a $k\times k $
Hankel matrix $M \in TP_2(c) $ with $\det M < 0 .$\\}

A simple consequence of Theorem 4 is the following fact

{\bf Corollary of Theorem 4.}{ \it   \\
(i) For every $1 \leq c < 4$ there exists a Toeplitz matrix $M \in
TP_2(c)$ but $M \notin TP; $ \\
(ii) for every $1 \leq c < 4$ there exists a Hankel matrix $M \in
TP_2(c)$ but $M \notin TP.$}

The following theorem shows that Theorem 1 remains valid for some
special classes of matrices with nonnegative elements.

{\bf Theorem 5.} {\it  Let  $M = (a_{i,j})$ be a  $k \times k$
matrix. Suppose that $\exists s,l \in {\bf Z}: -(k-1)\leq s <l\leq
k-1$ such that $a_{i,j}>0$ for $s\leq j-i \leq l$ and $a_{i,j}=0$
for $j-i<s$ or $j-i>l$ . If $a_{i,j}a_{i+1,j+1} \geq c_k \
a_{i,j+1}a_{i+1,j} \quad (1 \leq i < m,\ 1 \leq j < n) $ then
$\det M \geq 0.$ }

We will show how to prove Theorem 5 in the section "Proof of
Theorem 4".

A variation of Theorem 3 for the class of Toeplitz matrices was
proved by J.~I.~Hutchinson in \cite{hut}. To formulate his
result we need some notions.

The class of $m$-times positive sequences consists of the
sequences $\{a_k\}_{k=0}^\infty$ such that all minors of the
infinite matrix

\begin {equation}
\label{mat}
 \left\|
  \begin{array}{ccccc}
   a_0 & a_1 & a_2 & a_3 &\ldots \\
   0   & a_0 & a_1 & a_2 &\ldots \\
   0   &  0  & a_0 & a_1 &\ldots \\
   0   &  0  &  0  & a_0 &\ldots \\
   \vdots&\vdots&\vdots&\vdots&\ddots
  \end{array}
 \right\|
\end {equation}
of order not greater than $m$ are non-negative. The class of
$m$-times positive sequences is denoted by $PF_m.$ A sequence is
called a multiply positive sequence if it is $m$-times positive
for some $m\in{\bf N}.$  A sequence
$\{a_k\}_{k=0}^\infty$ such that all minors of the infinite matrix
(\ref{mat}) are nonnegative is called a totally positive sequence.
The class of totally positive sequences is denoted by $PF_\infty.$
The corresponding classes of generating functions
$$f(z)=\sum_{k=0}^\infty a_kz^k$$
are also denoted by $PF_m$ and $PF_\infty$.

The multiply positive sequences (also called P\'olya frequency
sequences) were introduced by Fekete in 1912 see~\cite{fek} in
connection with the problem of exact calculation of the number of
positive zeros of a real polynomial.

The class $PF_\infty$ was completely described by Aissen,
Schoenberg, Whitney and Edrei in~\cite{aissen} (see also
\cite[p.412]{tp}):

{\bf Theorem ASWE. } {\it A function $f \in PF_\infty$ iff
$$f(z)=C z^n e^{\gamma z}\prod_{k=1}^\infty (1+\alpha_kz)/(1-\beta_kz),$$
where $C\ge 0, n\in {\bf Z},
\gamma\ge0,\alpha_k\ge0,\beta_k\ge0,\sum(\alpha_k+\beta_k)
<\infty.$}

By Theorem ASWE a polynomial $p(z) = \sum _{k=0}^n a_k z^k, \  a_k
\geq 0 ,$ has only real zeros if and only if the sequence $(a_0,
a_1, \ldots , a_n, 0, 0, \ldots ) \in PF_\infty .$

In 1926, Hutchinson \cite[p.327]{hut} extended the work of
Petrovitch \cite{pet} and Hardy \cite{har1} or \cite[pp. 95-100]{har2}
and proved the following theorem.

{\bf Theorem B.} { \it Let $f(z)=\sum_{k=0}^\infty a_k z^k$,
$a_k>0,\ \forall k.$ Inequality
\begin{equation}
\label{ner}a_n^2\geq 4 a_{n-1}a_{n+1},\
\forall n\geq 1
\end{equation}
holds if and only if the following two properties hold:\\
(i) The zeros of f(x)
are all real, simple and negative and \\
(ii) the zeros of any polynomial $\sum_{k=m}^n a_kz^k$, formed by
taking any number of consecutive terms of $f (x)$, are all real
and non-positive.}

It is easy to see that (\ref{ner}) implies
$$a_n \leq \frac
{a_1}{4^{n(n-1)/2}} \left(
\frac{a_1}{a_0} \right) ^{n-1},\  n\geq 2,$$ that is $f$ is an
entire function of the order $0.$ So by the Hadamard theorem
(see, for example, \cite[p. 24]{lev})
$$f(z) =  C z^n \prod_{k=1}^\infty (1+\alpha_kz),$$
where $C\ge 0, n\in {\bf N \cup \{ 0\}},
\alpha_k\geq 0,\sum(\alpha_k)
<\infty. $

Using ASWE Theorem we obtain from Theorem B that
\begin{equation}
\label{123} a_n^2\geq 4 \  a_{n-1}a_{n+1},\ \forall n\geq 1
\Rightarrow \{a_n \}_{n=0}^\infty \in PF_\infty.
\end{equation}

In \cite{klv} it was proved that the constant $4$ in (\ref{123})
is sharp.

Thus, Theorem B provides a simple sufficient condition for deducing
when a sequence is a totally positive sequence. Theorem 5 provides
the following simple
sufficient condition of multiple positivity for a sequence.

{\bf Corollary of Theorem 5.}  {\it  Let $\{a_n \}_{n=0}^\infty$
 be a sequence of nonnegative numbers. Then
 $$a_n^2\geq c_m a_{n-1}a_{n+1},\
\forall n\geq 1   \Rightarrow \{a_n \}_{n=0}^\infty \in PF_m. $$}

Our results are applicable also to the moment problem.
Recall that a sequence of positive
numbers $\{s_k\}_{k=0}^\infty $ is said to be the moment sequence
of a nondecreasing function $F :
{\bf R} \rightarrow {\bf R}$ if
$$ s_k = \int_{-\infty}^\infty t^k  d F(t).$$
A sequence of positive numbers is called a Hamburger moment
sequence if it is a moment sequence of a function $F$
having infinitely many points of growth. The following famous
theorem gives the description of Hamburger moment sequences.

{\bf Theorem C.} (\cite{hamb}, see also \cite[chapt. 2]{akh})
{ \it A sequence of positive
numbers $\{s_k\}_{k=0}^\infty $ is a Hamburger moment
sequence if and only if
\begin {equation}
 \det \left(
  \begin{array}{cccc}
   s_0 & s_1 & \ldots  &s_k \\
   s_1   & s_2 &\ldots & s_{k+1} \\
   \vdots &  \vdots  & \ldots & \vdots \\
   s_k  &  s_{k+1}  &  \ldots  & s_{2k} \\
    \end{array}
 \right)> 0, \quad k=0, 1, 2, \ldots .
\end {equation} }

The following statement is proved in \cite{bis}.

{\bf Theorem D.} {\it Let $ d$ be the positive
solution of $\sum_{n=1}^\infty d^{- n^2}=1/4 $ $(d \approx 4{.}06) .$
Then any positive sequence $\{s_k\}_{k=0}^\infty $ satisfying
$$s_{n-1}s_{n+1}\geq d s_n^2 \quad n = 0, 1, 2, \ldots $$
is a Hamburger moment sequence.}

Theorem 3 implies the following statement.

{\bf Corollary of Theorem 3.}  {\it Any positive sequence
$\{s_k\}_{k=0}^\infty $ satisfying
$$s_{n-1}s_{n+1}\geq 4 s_n^2 \quad n = 0, 1, 2, \ldots $$
is a Hamburger moment sequence.}

The constant $4$ in the Corollary above cannot be improved.

\section{Proof of Theorem 1.}

We need the following sequence of functions:
\begin{equation}
\label{d0}
F_m (c) = \sum_{j=0}^{\lfloor m/2\rfloor} {m-j\choose j} (-1)^j \frac{1}{c^j},
\quad m=0,1, 2, \ldots , \ c\geq 1,
\end{equation}
where by $\lfloor x \rfloor$ we denote the integral part of $x.$

The following lemma provides  some properties for this sequence of functions.

{\bf Lemma 1.} {\it  \\
(i) The following identities hold
\begin {equation}
\label{d1}
\begin{array}{cc}
& F_0(c) = F_1(c) = 1 \\
\nonumber & F_m(c) = F_{m-1}(c) - \frac{1}{c} F_{m-2}(c),
\quad m=2, 3, 4,  \ldots .
\end{array}
\end{equation}

(ii) For $c = 4 \cos^2 \phi $ we have
\begin{equation}
\label{d2}
F_m (c) = \frac {\sin(m+1)\phi}{ c^{m/2}  \sin\phi}.
\end{equation} }

(iii) For $c_k = 4 \cos^2 \frac{\pi}{k+1} $ we have
\begin{equation} \label{t4}
F_{j-1}(c_k)-\frac{1}{c_k^2}F_{j-2}(c_k)-\frac{1}{c_k^j} \geq
F_j(c_k),\quad k\geq 3, \quad j=2,3,\ldots, k-1.
\end{equation}

{\it Proof of Lemma 1.} Formula (\ref{d1}) follows directly from
(\ref{d0}). Formula (\ref{d2}) is a simple consequence of the
well-known trigonometric identity (see, for example, \cite[p. 696]{zwill})
$$ \frac {\sin(m+1)\phi}{  \sin\phi} = \sum_{j=0}^{\lfloor m/2\rfloor}
 {m-j\choose j} (-1)^j
(2\cos\phi)^{m-2j}.$$

Using the identity $4\cos^2\phi-1 =\frac{\sin(3\phi)}{\sin\phi}$
we have
$$ F_{j-1}(c_k)-\frac{1}{c_k^2}F_{j-2}(c_k)-\frac{1}{c_k^j} -
F_j(c_k)=(\frac{1}{c_k}-\frac{1}{c_k^2})F_{j-2}(c_k)-
\frac{1}{c_k^j}$$
$$ =\frac{1}{c_k^{(j+2)/2}}\left ( \frac {\sin
(3\frac{\pi}{k+1})}{\sin\frac{\pi}{k+1}}\cdot \frac {\sin
((j-1)\frac{\pi}{k+1})}{\sin\frac{\pi}{k+1}}-\frac{1}
{(2\cos\frac{\pi}{k+1})^{j-2}}\right )\geq$$
$$\frac{1}{c_k^{(j+2)/2}}\left ( \frac {\sin
(3\frac{\pi}{k+1})}{\sin\frac{\pi}{k+1}}\cdot \frac {\sin
((j-1)\frac{\pi}{k+1})}{\sin\frac{\pi}{k+1}}-1\right )\geq 0,$$
for $k\geq 3$ and $j=2,3, \ldots, k-1.$ Inequality (\ref{t4}) is
proved.

Lemma 1 is proved. $\Box$

The following Lemma was proved in \cite{cc}.

{\bf Lemma A.} {\it  Let  $M = (a_{i,j}), \ 1\leq i \leq m,\ 1 \leq
j \leq n $ and $M\in TP_2(c) , c\geq 1. $ Then
$$ a_{i,j}a_{k,l} \geq c^{(l-j)(k-i)} a_{i,l}a_{k,j} ,
\quad \mbox{for all} \quad i<k,  j<l. $$}

A simple consequence of Lemma A is the fact that if $M\in TP_2(c)$
then any submatrix of $M$ also belongs to $TP_2(c).$ Analogously
if $M\in STP_2(c)$ then any submatrix of $M$ also belongs to
$STP_2(c).$

For a matrix $M =(a_{i,j})$ we will denote by
$M\left(^{i_1,i_2,\ldots, i_k}_{j_1,j_2,\ldots, j_k}\right )$
the following submatrix of $M$
$$ M\left(^{i_1,i_2,\ldots, i_k}_{j_1,j_2,\ldots, j_k}\right ) =
\left( \begin{array}{cccc}
   a_{i_1,j_1} & a_{i_1,j_2} & \ldots & a_{i_1,j_k} \\
   a_{i_2,j_1} & a_{i_2,j_2} & \ldots & a_{i_2,j_k} \\
  \vdots & \vdots  & \cdots & \vdots\\
   a_{i_k,j_1} & a_{i_k,j_2} & \ldots & a_{i_k,j_k}\\
   \end{array} \right )$$

We now prove the following claim (which consists of three parts)
by induction on $n.$ Let  $M = (a_{i,j})$ be an  $n
\times n$ matrix and $M\in TP_2(c),$ where $c\geq 4\cos ^2
\frac{\pi}{n+1}.$ Then the following inequalities hold:
\begin{equation}
\label{h1}
\det M \geq 0 .
\end{equation}
\begin{equation}
\label{h2}
\det M \geq a_{1,1}\det M\left(^{2,3,\ldots, n}_{2,3,\ldots, n}\right )
- a_{1,2}a_{2,1}\det M\left(^{3,4, \ldots, n}_{3,4, \ldots, n}\right ).
\end{equation}
\begin{equation}
\label{h3}
\det M \leq a_{1,1}\det M\left(^{2,3,\ldots, n}_{2,3,\ldots, n}\right ).
\end{equation}

Since $M \in TP_2(c)$ then hypothesis (\ref{h1}), (\ref{h2}),
(\ref{h3}) are true for $n=2.$ The proof below is based on the
following lemma.

{\bf Lemma 2.} { \it Let $c_0 \geq 1,$ $M=(a_{i,j})\in TP_2(c_0)$
be an $n\times n$ matrix satisfying the following conditions \\
(i) $\forall i=2, 3, \ldots, n \quad \det M\left(^{i, i+1,\ldots,
n}_{i,i+1,\ldots, n}\right ) \geq 0 ;$ \\
(ii) $\forall i=1, 2, \ldots, n-2 $
$$\det M\left(^{i, i+1,\ldots, n}_{i, i+1,\ldots, n}\right )\geq
a_{i,i}\det M\left(^{i+1, i+2,\ldots, n}_{i+1, i+2,\ldots, n}\right
)-a_{i,i+1}a_{i+1,i} \det M\left(^{i+2, i+3,\ldots, n}_{i+2, i+3,
\ldots, n}\right ).$$
Then for all $c, 1\leq c \leq c_0$ the
following inequalities are valid:}
\begin{eqnarray}
\label{l1}
&\det M\left(^{m+1, m+2, \ldots, n}_{m+1, m+2, \ldots,
n}\right )\geq a_{m+1, m+1} \left (\det M\left(^{m+2, m+3, \ldots, n}_{m+2, m+3,
\ldots, n}\right ) - \right. \\
\nonumber & \left. \frac {1}{c}a_{m+2,m+2}\det M\left(^{m+3,
m+4, \ldots, n}_{m+3, m+4, \ldots, n}\right )\right ),\   m=0, 1,
\ldots, n-3.
\end{eqnarray}

\begin{eqnarray}
\label{l2}
&\det M\geq a_{1,1}a_{2,2} \cdots  a_{m, m}
\left (F_m(c)\det M\left(^{m+1, m+2, \ldots, n}_{m+1,
m+2, \ldots, n}\right ) - \right. \\
\nonumber & \left. \frac {1}{c}F_{m-1}(c) a_{m+1,m+1}\det
M\left(^{m+2, m+3, \ldots, n}_{m+2, m+3, \ldots, n}\right )\right
), m= 1, 2, \ldots, n-2.
\end{eqnarray}

\begin{eqnarray}
\label{l3}
&F_m(c)\det M\left(^{m+1, m+2, \ldots, n}_{m+1, m+2,
\ldots, n}\right ) - \frac {1}{c} F_{m-1}(c) a_{m+1,m+1}\det
M\left(^{m+2, m+3, \ldots, n}_{m+2, m+3, \ldots, n}\right )\geq \\
\nonumber & a_{m+1,m+1}\left (F_{m+1}(c)\det M\left(^{m+2, m+3,
\ldots, n}_{m+2, m+3, \ldots, n}\right ) - \frac {1}{c}F_m(c)
a_{m+2,m+2}\det M\left(^{m+3, m+4, \ldots, n}_{m+3, m+4, \ldots,
n}\right )\right ),\\
\nonumber & m=1, 2, \ldots, n-3.
\end{eqnarray}

\begin{eqnarray}
\label{l4} &F_m(c)\det M\left(^{m+1, m+2, \ldots, n}_{m+1, m+2,
\ldots, n}\right ) - \frac {1}{c}F_{m-1}(c) a_{m+1,m+1}\det
M\left(^{m+2, m+3, \ldots, n}_{m+2, m+3, \ldots, n}\right )\geq
\\
\nonumber &a_{m+1,m+1}a_{m+2,m+2} \cdots a_{n,n}
F_n(c),\quad m=1, 2, \ldots, n-2.
\end{eqnarray}

{\it Proof of Lemma 2.} First we prove (\ref {l1}). Since
$M\in TP_2(c)$ and by (ii) we have
$$ \det M\left(^{m+1, m+2, \ldots, n}_{m+1, m+2, \ldots,
n}\right )\geq a_{m+1, m+1} \det M\left(^{m+2, m+3, \ldots,
n}_{m+2, m+3, \ldots, n}\right ) - $$
$$a_{m+1,m+2}a_{m+2,m+1}\det
M\left(^{m+3, m+4, \ldots, n}_{m+3, m+4, \ldots, n}\right )\geq
a_{m+1, m+1} \det M\left(^{m+2, m+3, \ldots, n}_{m+2, m+3, \ldots,
n}\right )$$
$$- \frac {1}{c}a_{m+1, m+1}a_{m+2,m+2}\det M\left(^{m+3, m+4,
\ldots, n}_{m+3, m+4, \ldots, n}\right ),\   m=0, 1, \ldots,
n-3.$$ Inequality (\ref {l1}) is proved.

Let us prove (\ref {l3}). Multiplying (\ref {l1}) by $F_m(c)$ we have
$$F_m(c)\det M\left(^{m+1, m+2, \ldots, n}_{m+1, m+2,
\ldots, n}\right ) - \frac {1}{c} F_{m-1}(c) a_{m+1,m+1}\det
M\left(^{m+2, m+3, \ldots, n}_{m+2, m+3, \ldots, n}\right )\geq
 $$
$$
a_{m+1,m+1} \left (\left(F_m(c)-\frac{1}{c}F_{m-1}(c)\right)\det M\left(^{m+2,
m+3, \ldots, n}_{m+2, m+3, \ldots, n}\right ) - \right.$$
$$\left. \frac {1}{c}F_m(c)
a_{m+2,m+2}\det M\left(^{m+3, m+4, \ldots, n}_{m+3, m+4, \ldots,
n}\right )\right ),  m=1, 2, \ldots, n-3,$$
and, using (\ref{d1})we obtain (\ref{l3}).

To prove (\ref{l4}) we apply (\ref{l3}) $(n-2-m)$ times. We derive
$$F_m(c)\det M\left(^{m+1, m+2, \ldots, n}_{m+1, m+2,
\ldots, n}\right ) - \frac {1}{c} F_{m-1}(c) a_{m+1,m+1}\det
M\left(^{m+2, m+3, \ldots, n}_{m+2, m+3, \ldots, n}\right )\geq$$
$$ a_{m+1,m+1}a_{m+2,m+2} \cdots  a_{n-2,n-2}\left
(F_{n-2}(c)\det M\left(^{n-1,n}_{n-1,n}\right ) - \frac
{1}{c}F_{n-3}(c) a_{n-1,n-1}a_{n,n}\right ).$$

Since $M\in TP_2(c_0)$ the following inequality holds for all $c,
1\leq c \leq c_0, $
\begin{equation}
\label{l5}
\\det M\left(^{n-1,n}_{n-1,n}\right )
\geq (1-\frac{1}{c})a_{n-1,n-1}a_{n,n},
\end{equation}
so by (\ref{d1}) we obtain
$$F_m(c)\det M\left(^{m+1, m+2, \ldots, n}_{m+1, m+2,
\ldots, n}\right ) - \frac {1}{c} F_{m-1}(c) a_{m+1,m+1}\det
M\left(^{m+2, m+3, \ldots, n}_{m+2, m+3, \ldots, n}\right )\geq$$
$$a_{m+1,m+1}a_{m+2,m+2} \cdots  a_{n,n} \left ( \left (
F_{n-2}(c)- \frac{1}{c}F_{n-3}(c) \right ) -
\frac{1}{c}F_{n-2}(c)\right ) =$$
$$ a_{m+1,m+1}a_{m+2,m+2} \cdots a_{n,n}
\left ( F_{n-1}(c) - \frac{1}{c}F_{n-2}(c)\right ) = $$
$$a_{m+1,m+1}a_{m+2,m+2}
\cdots a_{n,n} F_n(c). $$ Inequality (\ref{l4}) is proved.

By (\ref{d1}) we rewrite inequality (\ref{l1}) for $m=0$ in the
following form:
$$\det M\geq a_{1,1}
\left (F_1(c)\det M\left(^{2, 3, \ldots, n}_{2, 3, \ldots,
n}\right ) - \frac {1}{c}F_0(c) a_{2,2}\det M\left(^{3, 4, \ldots,
n}_{3, 4, \ldots, n}\right )\right ).$$

To prove (\ref{l2}) we apply (\ref{l3}) $(m-1)$ times.

Lemma 2 is proved. $\Box$

{\bf Remark.} If a matrix $M$ satisfies the conditions of Lemma 2
and, moreover,  $a_{n-1,n-1}a_{n,n} > c_0 a_{n-1,n}a_{n,n-1},$
then inequality (\ref{l5}) is strict, hence (\ref{l4}) is strict,
i.e.
\begin{eqnarray}
\label{l6} &F_m(c)\det M\left(^{m+1, m+2, \ldots, n}_{m+1, m+2,
\ldots, n}\right ) - \frac {1}{c}F_{m-1}(c) a_{m+1,m+1}\det
M\left(^{m+2, m+3, \ldots, n}_{m+2, m+3, \ldots, n}\right )>
\\
\nonumber &a_{m+1,m+1}a_{m+2,m+2} \cdots  a_{n,n}
F_n(c),\quad m=1, 2, \ldots, n-2.
\end{eqnarray}
In particular, for all matrices $M\in STP(c_0)$ inequality
(\ref{l6}) is valid for all  $c, 1\leq c \leq c_0.$

Assume that conditions (\ref{h1}), (\ref{h2}) and (\ref{h3}) hold
for all matrices of sizes smaller than $k.$ Let us prove these
conditions for $n=k.$

{\bf Lemma 3. } {\it  Let $M=(a_{i,j})$ be a $k\times k$ matrix,
$M\in TP_2(c)$, $c\geq c_k := 4\cos ^2 \frac{\pi}{k+1}.$
For all $j=2,3, \ldots, k-1$ the following
inequality holds. }
$$ a_{1,j}\det M\left(^{2, 3, \ldots, k}_{1, 2, \ldots,j-1,j+1,
\ldots, k}\right )- a_{1,j+1}\det M\left(^{2, 3, \ldots, k}_{1, 2,
\ldots,j,j+2, \ldots, k}\right )\geq 0. $$

{\it Proof of Lemma 3.} Since $m\in TP_2(c),  $ $ M\left(^{2, 3,
\ldots, k}_{1, 2, \ldots,j-1,j+1, \ldots, k}\right )\in TP_2(c) $
and $ M\left(^{2, 3, \ldots, k}_{1, 2, \ldots,j,j+2, \ldots,
k}\right )\in TP_2(c).$ Since $4\cos ^2 \frac{\pi}{n+1} \leq 4\cos
^2 \frac{\pi}{k+1}$ for $n=2,3, \ldots, k-1$ we can apply the
induction hypothesis to the matrices $M\left(^{2, 3, \ldots,
k}_{1, 2, \ldots,j-1,j+1, \ldots, k}\right ), M\left(^{2, 3,
\ldots, k}_{1, 2, \ldots,j,j+2, \ldots, k}\right )$ and to all
their square submatrices. We apply inequality (\ref{h3}) $j$ times
and obtain
$$\det M\left(^{2, 3, \ldots, k}_{1, 2, \ldots,j,j+2, \ldots,
k}\right )\leq a_{2,1}a_{3,2} \cdots  a_{j+1,j}\det
M\left(^{j+2, j+3, \ldots, k}_{j+2, j+3, \ldots,k}\right ).$$
From Lemma A and from the fact
$a_{1,j+1}a_{j+1,j} \leq  \frac{1}{c_k^j} a_{1,j} a_{j+1,j+1} $
now we conclude
\begin{equation}
\label{t1} a_{1,j+1}\det M\left(^{2, 3, \ldots, k}_{1, 2,
\ldots,j,j+2, \ldots, k}\right )\leq
\frac{1}{c_k^j}a_{1,j}a_{2,1}a_{3,2} \cdots
a_{j,j-1}a_{j+1,j+1}\det M\left(^{j+2, j+3, \ldots, k}_{j+2, j+3,
\ldots,k}\right ).
\end{equation}

By the induction hypothesis the matrix $M\left(^{2, 3, \ldots,
k}_{1, 2, \ldots, j-1,j+1, \ldots , k}\right )$ satisfies the
assumptions of Lemma 2. Applying to this matrix (\ref{l2}) with
$m=j-2$ we obtain
$$\det M\left(^{2, 3, \ldots, k}_{1, 2, \ldots,
j-1,j+1, \ldots , k}\right )\geq a_{2,1}a_{3,2} \cdots
a_{j-1, j-2} \left (F_{j-2}(c_k)\det M\left(^{j,
j+1,j+2 \ldots, k}_{j-1, j+1,j+2 \ldots, k}\right )\right. $$
$$ \left. - \frac
{1}{c_k}F_{j-3}(c_k) a_{j,j-1}\det M\left(^{j+1, j+2, \ldots,
k}_{j+1, j+2, \ldots, k}\right )\right ).$$
Applying (\ref{h2}) to the matrix $M\left(^{j, j+1,j+2 \ldots,
k}_{j-1, j+1,j+2 \ldots, k}\right )$ and plugging the result into the last
formula we have
$$\det M\left(^{2, 3, \ldots, k}_{1, 2, \ldots,
j-1,j+1, \ldots , k}\right )\geq a_{2,1}a_{3,2} \cdots
a_{j-1, j-2}  \left( a_{j,j-1}(F_{j-2}(c_k)-\right. $$
$$\left. \frac{1}{c_k}F_{j-3}(c_k))
\det M\left(^{j+1,j+2 \ldots, k}_{j+1,j+2 \ldots,
k}\right )- a_{j,j+1}a_{j+1,j-1}F_{j-2}(c_k) \det M\left(^{j+2,
j+3, \ldots, k}_{j+2, j+3, \ldots, k}\right ) \right ),$$
whence, by Lemma A and (\ref{d1})  we obtain
$$\det M\left(^{2, 3, \ldots, k}_{1, 2, \ldots,
j-1,j+1, \ldots , k}\right )\geq a_{2,1}a_{3,2} \cdots
a_{j-1, j-2}a_{j,j-1}
\left ( F_{j-1}(c_k)\det M\left(^{j+1,j+2 \ldots, k}_{j+1,j+2 \ldots,
k}\right )- \right. $$
$$\left. \frac{1}{c_k^2} a_{j+1,j+1}F_{j-2}(c_k) \det
M\left(^{j+2, j+3, \ldots, k}_{j+2, j+3, \ldots, k}\right )\right
).$$
Further applying (\ref{l1}) to the matrix $M\left(^{j+1,j+2
\ldots, k}_{j+1,j+2 \ldots, k}\right )$ we have
\begin{eqnarray}
\label{t2} &\det M\left(^{2, 3, \ldots, k}_{1, 2, \ldots, j-1,j+1,
\ldots , k}\right )\geq a_{2,1}a_{3,2} \cdots
a_{j,j-1}a_{j+1,j+1} \left(  \det M\left(^{j+2,j+3 \ldots,
k}_{j+2,j+3 \ldots, k}\right )\right.
\\
\nonumber &
\left.(F_{j-1}(c_k)-\frac{1}{c_k^2}F_{j-2}(c_k))- \frac{1}{c_k}
a_{j+2,j+2}F_{j-1}(c_k) \det M\left(^{j+3, j+4, \ldots, k}_{j+3,
j+4, \ldots, k}\right ) \right).
\end{eqnarray}
By (\ref{t1}) and (\ref{t2}) we derive
\begin{eqnarray}
\label{t3} &a_{1,j}\det M\left(^{2, 3, \ldots, k}_{1, 2, \ldots,
j-1,j+1, \ldots , k}\right ) - a_{1,j+1}\det M\left(^{2, 3,
\ldots, k}_{1, 2, \ldots, j,j+2, \ldots , k}\right )\\
\nonumber &\geq a_{1,j}a_{2,1}a_{3,2} \cdots
a_{j,j-1}a_{j+1,j+1}\left((
F_{j-1}(c_k)-\frac{1}{c_k^2}F_{j-2}(c_k)-\frac{1}{c_k^j})\right. \\
\nonumber & \left.  \det M\left(^{j+2,j+3 \ldots,
k}_{j+2,j+3 \ldots, k}\right )- \frac{1}{c_k}
a_{j+2,j+2}F_{j-1}(c_k) \det M\left(^{j+3, j+4, \ldots, k}_{j+3,
j+4, \ldots, k}\right ) \right).
\end{eqnarray}

It follows from (\ref{t3}),(\ref{t4}) and (\ref{l4}) that
$$a_{1,j}\det M\left(^{2, 3, \ldots, k}_{1, 2, \ldots, j-1,j+1,
\ldots , k}\right ) - a_{1,j+1}\det M\left(^{2, 3, \ldots, k}_{1,
2, \ldots, j,j+2, \ldots , k}\right )\geq$$
 $$a_{1,j}a_{2,1}a_{3,2} \cdots
a_{j,j-1}a_{j+1,j+1}\left(F_j(c_k) \det M\left(^{j+2,j+3 \ldots,
k}_{j+2,j+3 \ldots, k}\right )- \right. $$
$$\frac{1}{c_k}
a_{j+2,j+2}F_{j-1}(c_k)\left. \det M\left(^{j+3, j+4,
\ldots, k}_{j+3, j+4, \ldots, k}\right ) \right )\geq $$
$$ a_{1,j}a_{2,1}a_{3,2} \cdots
a_{j,j-1}a_{j+1,j+1}a_{j+2,j+2} \cdots
a_{k,k}F_{k-1}(c_k).
$$
Hence by Lemma 1 and (\ref{d2}) with $m=k-1$ we conclude that
$$a_{1,j}\det M\left(^{2, 3, \ldots, k}_{1, 2, \ldots, j-1,j+1,
\ldots , k}\right ) - a_{1,j+1}\det M\left(^{2, 3, \ldots, k}_{1,
2, \ldots, j,j+2, \ldots , k}\right )$$ $$\geq
a_{1,j}a_{2,1}a_{3,2} \cdots
a_{j,j-1}a_{j+1,j+1}a_{j+2,j+2} \cdots  a_{k,k} \frac
{\sin (k\frac{\pi}{k+1})}{c_k^{(k-1)/2}\sin\frac{\pi}{k+1}} \geq
0.$$

Lemma 3 is proved. $\Box$

Now we will prove (\ref{h2}). Using Lemma 3 we have
$$\det M\left(^{1, 2, \ldots, k}_{1, 2, \ldots,  k}\right )
= \sum_{j=1}^k (-1)^{j+1} a_{1,j}\det M\left(^{2, 3, \ldots, k}_{1, 2,
\ldots, j-1,j+1, \ldots , k}\right )$$  $$ \geq a_{1,1}\det
M\left(^{2, 3, \ldots, k}_{ 2,3 \ldots, k}\right ) - a_{1,2}\det
M\left(^{2, 3, \ldots, k}_{1, 3,4,  \ldots, k}\right ). $$ We
apply the induction hypothesis (\ref{h3}) to the matrix
$M\left(^{2, 3, \ldots, k}_{1, 3,4,  \ldots, k}\right ).$ We have
$$\det M\left(^{1, 2, \ldots, k}_{1, 2, \ldots,  k}\right )
\geq a_{1,1}\det M\left(^{2, 3, \ldots, k}_{ 2,3 \ldots, k}\right
) - a_{1,2}a_{2,1}\det M\left(^{3, 4, \ldots, k}_{ 3,4,  \ldots,
k}\right ). $$
The inequality (\ref{h2}) is proved.

By Lemma 3
$$\det M\left(^{1, 2, \ldots, k}_{1, 2, \ldots,  k}\right )
= \sum_{j=1}^k (-1)^{j+1} a_{1,j}\det M\left(^{2, 3, \ldots, k}_{1, 2,
\ldots, j-1,j+1, \ldots , k}\right ) \leq a_{1,1}\det M\left(^{2,
3, \ldots, k}_{ 2,3 \ldots, k}\right ) .$$
The inequality (\ref{h3}) is proved.

To prove (\ref{h1}) we note that by (\ref{h2}) and induction
hypothesis the matrix $M$ satisfies the assumptions of Lemma 2. It
follows from (\ref{l2}), (\ref{l4}) and Lemma 1 that
$$ \det M \geq a_{1,1}a_{2,2} \cdots  a_{k,k}F_k(c_k) =
a_{1,1}a_{2,2} \cdots  a_{k,k} \frac {\sin
\pi}{c_k^{k/2}\sin\frac{\pi}{k+1}}=0.$$
Hence the statement (i) in Theorem 1 is proved.

Now we will prove the statement (ii) in Theorem 4. If $M\in
STP_k(c_k)$ then by (\ref{l6}) we can rewrite the last inequality
in the following form
$$ \det M > a_{1,1}a_{2,2} \cdots  a_{k,k}F_k(c_k) =
a_{1,1}a_{2,2} \cdots  a_{k,k} \frac {\sin
\pi}{c_k^{k/2}\sin\frac{\pi}{k+1}}=0.$$
Hence the statement (ii) in Theorem 1 is proved, which completes the proof
of Theorem 1. $\Box$

In fact, we have proved a slightly stronger theorem, which may be of
independent interest.

{\bf Theorem 6.} {\it Suppose $c \geq 4\cos^2\frac{\pi}{k+1}.$ Let
$M=(a_{i,j})\in TP_2(c)$ be a  $k\times k$ matrix.   Then
$$\det M \geq a_{1,1}a_{2,2} \cdots  a_{k,k}F_k(c).$$}

\section{Proof of Theorem 4.}

Note that $TP_2(c_1)\subset TP_2(c_2)$ for $c_1\geq c_2$. Thus it is
sufficient to prove Theorem~4 with $c \in (c_k -\varepsilon, c_k)$
for $\varepsilon > 0$ being small enough.

Consider the following $n\times n$ symmetrical Toeplitz matrix.
\begin {equation}
M_n(\phi):= \left\|
  \begin{array}{ccccccc}
   2\cos \phi & 1 & 0 & 0 &\ldots &0&0\\
   1   & 2\cos \phi & 1 & 0 &\ldots&0&0 \\
   0   &  1  & 2\cos \phi & 1 &0&\ldots&0 \\
   \vdots&\vdots&\vdots&\vdots&\ldots&\vdots&\vdots\\
   0&0&\ldots&0&1   & 2\cos \phi & 1\\
   0  &  0  &0 &\ldots&0& 1 & 2\cos \phi \\
     \end{array}
 \right\|,
\label{g1}
\end {equation}

where $0 \leq \phi < \pi/2.$ Obviously, $M_n(\phi)\in
TP_2(4\cos^2\phi).$ The matrix $M_n(\phi)$ satisfies the following
recursion relation $\det M_n(\phi) = 2\cos\phi \det
M_{n-1}(\phi) - \det M_{n-2}(\phi) $ and $M_1(\phi) =2\cos\phi,
M_2(\phi)= 4\cos^2\phi -1.$ It is easy to verify that $\det
M_n(\phi) = \frac{\sin(n+1)\phi}{\sin \phi}.$ So for all $\phi \in
(\frac{\pi}{n+1}, \frac{2\pi}{n+1})$ we have $\det M_n(\phi) <0. $
For $\phi \in (\frac{\pi}{n+1}, \frac{2\pi}{n+1})$ consider the
following $n\times n$ symmetrical  Toeplitz matrix
\begin {eqnarray}
& T_n(\phi , \varepsilon_1,\ldots , \varepsilon_{n-2}):= \\
\nonumber & \left\|
  \begin{array}{ccccccc}
   2\cos \phi & 1 & \varepsilon_1 & \varepsilon_2 &\ldots &\varepsilon_{n-3}
   &\varepsilon_{n-2}\\
   1   & 2\cos \phi & 1 & \varepsilon_1 &\ldots&
   \varepsilon_{n-4}&\varepsilon_{n-3} \\
   \varepsilon_1   &  1  & 2\cos \phi & 1 &\varepsilon_1&\ldots&
   \varepsilon_{n-4} \\
   \vdots&\vdots&\vdots&\vdots&\ldots&\vdots&\vdots\\
   \varepsilon_{n-3}&\varepsilon_{n-4}&\ldots&\varepsilon_1&1&
    2\cos \phi & 1\\
   \varepsilon_{n-2}  &  \varepsilon_{n-3}  &\varepsilon_{n-4} &
   \ldots&\varepsilon_1& 1 & 2\cos \phi \\
     \end{array}
 \right\|,
\label{g11}
\end {eqnarray}
where  $\varepsilon_1 > \varepsilon_2 > \cdots >
\varepsilon_{n-2}>0$ and $\varepsilon_1$ is chosen to satisfy the
inequality $1\geq 4\cos^2\phi \cdot 2\cos \phi \cdot
\varepsilon_1,$ then $\varepsilon_2$ is chosen to satisfy the
inequality $\varepsilon_1^2\geq 4\cos^2\phi  \cdot \varepsilon_2,
 $ then $\varepsilon_3$ is chosen to satisfy the inequality
$\varepsilon_2^2\geq 4\cos^2\phi  \cdot \varepsilon_1\cdot
\varepsilon_3,\ \ldots $ and then $\varepsilon_{n-2}$ is chosen to
satisfy the inequality $\varepsilon_{n-3}^2\geq 4\cos^2\phi  \cdot
\varepsilon_{n-4}\cdot \varepsilon_{n-2}.$ Under these conditions
we have $T_n(\phi , \varepsilon_1,\ldots , \varepsilon_{n-2}) \in
TP_2(4\cos^2\phi).$ Since $T_n(\phi , 0,0, \ldots , 0)=M_n(\phi)$
we obtain  $\det T_n(\phi , 0,0, \ldots , 0)< 0$ for $\phi \in
(\frac{\pi}{n+1}, \frac{2\pi}{n+1}).$ Therefore we have $\det
T_n(\phi , \varepsilon_1,\ldots , \varepsilon_{n-2})< 0$ for $\phi
\in (\frac{\pi}{n+1}, \frac{2\pi}{n+1})$ if $\varepsilon_1$ is
small enough.

Thus, for every $c\in (4\cos^2 \frac{2\pi}{n+1}, c_n)$ the
statement (i) of Theorem 4 is proved. Since $TP_2(c_1) \subset
TP_2(c_2)$ for $c_1 \geq c_2$ the statement (i) of Theorem 4
follows.

We use the same method to obtain the proof of Theorem 5.

To prove the statement (ii) we consider the following Hankel
matrix $D_n(p,q)$ with $p\geq 1, q\geq 1.$

\begin {equation}
D_n(p,q):= (
p^{\lfloor(i+j-2)/2\rfloor \lfloor(i+j-1)/2\rfloor}
q^{\lfloor(i+j-3)/2\rfloor \lfloor(i+j-2)/2\rfloor},\ 1\leq i, j
\leq n,
 \label{g2}
\end {equation}

or,

\begin {equation}
D_n(p,q) = \left\|
  \begin{array}{ccccccc}
   1 & 1 & p & p^2q &\ldots &\ast&\ast\\
   1 & p & p^2q & p^4 q^2 &\ldots&\ast&\ast \\
   p   &  p^2q  & p^4q ^2 & p^6q^4 &\ldots &\ast&\ast \\
   \vdots&\vdots&\vdots&\vdots&\ldots&\vdots&\vdots\\
   \ast&\ast&\ast&\ast&\ldots& p^{(n-2)^2}q^{(n-2)(n-3)} &
   p^{(n-1)(n-2)}q^{(n-2)^2}\\
   \ast&\ast&\ast&\ast&\ldots& p^{(n-1)(n-2)}q^{(n-2)^2} &
   p^{(n-1)^2}q^{(n-1)(n-2)} \\
     \end{array}
 \right\|,
\end {equation}

By direct calculation we obtain $D_n(p,q)\in TP_2(\min(p,q)).$

{\bf Lemma 4.}  {\it For all $n\geq 3$ we have
\begin {equation}
\label{ao1} \det D_n(p,q)= p^{\beta_n}q^{\alpha_n}F_n(p) +
 Q_{\alpha_n-1}(p,q),
\end {equation}
 where $\alpha_n=\frac{n(n-1)(n-2)}{3},  \beta_n = \frac{n(n-1)(2n-1)}{6}$
 and $Q_{\alpha_n-1}(p,q)$ is a polynomial in $p,q$ such that
 $\deg_q Q_{\alpha_n-1}(p,q) \leq \alpha_n-1 .$
 (Here and further by  $\deg_q Q(p,q)$ we will denote the degree of
  $Q(p,q)$ with respect to $q.$)}

{\it Proof.} We will prove this lemma by induction in $n.$ For
$n=3$ the statement is true as can be verified directly. The expansion of  $\det D_n(p,q)$
along column $n$ gives
\begin {eqnarray}
&\det D_n(p,q) =  R_{\alpha_n-1}(p,q) +\\
\nonumber
&\det \left\|
  \begin{array}{ccccccc}
   1 & 1 & p & p^2q &\ldots &\ast&0\\
   1 & p & p^2q & p^4 q^2 &\ldots&\ast&0 \\
   p   &  p^2q  & p^4q ^2 & p^6q^4 &\ldots &\ast&0\\
   \vdots&\vdots&\vdots&\vdots&\ldots&\vdots&\vdots\\
   \ast&\ast&\ast&\ast&\ldots& \ast &
   0\\
   \ast&\ast&\ast&\ast&\ldots& p^{(n-2)^2}q^{(n-2)(n-3)} &
   p^{(n-1)(n-2)}q^{(n-2)^2}\\
   \ast&\ast&\ast&\ast&\ldots& p^{(n-1)(n-2)}q^{(n-2)^2} &
   p^{(n-1)^2}q^{(n-1)(n-2)} \\
     \end{array}
 \right\| ,
\end {eqnarray}

where $R_{\alpha_n-1}(p,q)$ is a polynomial in $p,q$ and
$\deg_qR_{\alpha_n-1}(p,q)\leq \alpha_n-1.$

The expansion of the determinant on the right-hand side of the
last equation along row $n$ gives
\begin {eqnarray}
&\det D_n(p,q) =  S_{\alpha_n-1}(p,q) +\\
\nonumber &\det \left\|
  \begin{array}{ccccccc}
   1 & 1 & p & p^2q &\ldots &\ast&0\\
   1 & p & p^2q & p^4 q^2 &\ldots&\ast&0 \\
   p   &  p^2q  & p^4q ^2 & p^6q^4 &\ldots &\ast&0\\
   \vdots&\vdots&\vdots&\vdots&\ldots&\vdots&\vdots\\
   \ast&\ast&\ast&\ast&\ldots& \ast &
   0\\
   \ast&\ast&\ast&\ast&\ldots& p^{(n-2)^2}q^{(n-2)(n-3)} &
   p^{(n-1)(n-2)}q^{(n-2)^2}\\
   0&0&0&\ldots&0& p^{(n-1)(n-2)}q^{(n-2)^2} &
   p^{(n-1)^2}q^{(n-1)(n-2)} \\
     \end{array}
 \right\| ,
\end {eqnarray}

where $S_{\alpha_n-1}(p,q)$ is a polynomial in $p,q$ and $\deg_q
S_{\alpha_n-1}(p,q)\leq \alpha_n-1.$

The last equation provides the following recursion relation
$$D_n(p,q) =  p^{(n-1)^2}q^{(n-1)(n-2)} D_{n-1}(p,q) -
p^{2(n-1)(n-2)}q^{2(n-2)^2}D_{n-2}(p,q) + T_{\alpha_n -1}(p,q),
$$
where $T_{\alpha_n-1}(p,q)$ is a polynomial in $p,q$ and $\deg_q
T_{\alpha_n-1}(p,q)\leq \alpha_n-1.$

Using the induction hypothesis and formula (\ref{d1})
we obtain the statement of Lemma 4.

Lemma 4 is proved. $\Box$

Note that $p^{\lfloor n/2\rfloor}F_n(p)$ is a polynomial in $p$ of degree
$\lfloor n/2 \rfloor.$ By (\ref{d2}) it has the following $\lfloor n/2\rfloor$ roots:
$$ 4\cos^2\frac{\pi}{n+1}, 4\cos^2\frac{2\pi}{n+1}, \ldots ,
4\cos^2\frac{\lfloor n/2 \rfloor \pi}{n+1}.$$
Obviously, $4\cos^2\frac{\pi}{n+1}$ is the largest root of this
polynomial. Hence for $p\in (4\cos^2\frac{2\pi}{n+1},
4\cos^2\frac{\pi}{n+1})$ we have $F_n(p) < 0.$

Let us fix an arbitrary $p_0\in (4\cos^2\frac{2\pi}{n+1},
4\cos^2\frac{\pi}{n+1}).$ Since
$$ \det D_n(p_0,q)= q^{\alpha_n}(p_0^{\beta_n}F_n(p_0) +
 q^{-\alpha_n}Q_{\alpha_n-1}(p_0,q)),$$
 where $Q_{\alpha_n-1}(p_0,q)$ is a polynomial in $q$ and
 $\deg Q_{\alpha_n-1}(p_0,q) \leq \alpha_n -1,$ for $q$ being
 large enough (and $q > p_0$) we obtain $D_n(p_0,q)\in TP_2(p_0)$
  but $\det D_n(p_0,q)<0.$

 Thus, for every $p\in (4\cos^2 \frac{2\pi}{n+1}, c_n)$ the
statement (ii) of Theorem 4 is proved. Since $TP_2(c_1) \subset
TP_2(c_2)$ for $c_1 \geq c_2$ the statement (ii) of Theorem 4
follows.

Theorem 4 is proved. $\Box$

{\bf Remark.} This is a revised version of the paper originally
submitted to the journal "Linear Algebra and its Applications" in
summer of 2004. Recently in the paper \cite{dim} the authors
formulated a conjecture which coincides with the statement proved
in our Theorem 1.

{\bf ACKNOWLEDGEMENT.} The authors are deeply grateful to
Professor V.M.~Kadets for valuable suggestions.

\end{document}